\RequirePackage{amsmath}
\documentclass[a4paper,12pt]{article}%

\usepackage{mathptmx}
\usepackage[english]{babel}

\usepackage{amsfonts}
\usepackage{amssymb}
\usepackage{amsthm}
\usepackage{graphicx}
\usepackage{pdfpages}
\usepackage{relsize}

\newcommand{\setx}[1]{ \{ #1 \} }
\newcommand{\setminx}[1]{ \setminus \{ #1 \} }

\newcommand{\mytext}[1]{ \: \textrm{#1} \: }
\newcommand{\mysetdescr}[2]{\left\{ #1 \: \left| \: #2 \right. \right\} }
\newcommand{\darr}{{\downarrow \,}}
\newcommand{\uarr}{{\uparrow \,}}
\newcommand{\ouarr}{{\uparrow_{_{_{\!\!\!\!\!\!\circ}}}}}\newcommand{\odarr}{{\downarrow^{^{\!\!\!\!\!\!\circ}}}}

\newcommand{\myN}{\mathbb{N}}

\newcommand{\myNkz}[1]{\setx{ 0, \ldots , #1 }}
\newcommand{\nz}{\setx{0,1,2}}

\def\rarr{\rightarrow}

\newcommand{\mf}[1]{\mathfrak{ #1 }}

\def\fN{\mf{N}}

\def\fNS{\mf{N}^{Seg}}
\def\fNL{\mf{N}^{LS}}
\def\fNU{\mf{N}^{US}}

\def\BP{\begin{proof}}
\def\EP{\end{proof}}

\DeclareMathOperator{\id}{id}

\begin{document}

\theoremstyle{plain}
\newtheorem{theorem}{Theorem}[section]
\newtheorem{definition}[theorem]{Definition}
\newtheorem{corollary}[theorem]{Corollary}
\newtheorem{lemma}[theorem]{Lemma}
\newtheorem{proposition}[theorem]{Proposition}
\newtheorem{criteria}[theorem]{Criteria}



\title{A contribution to the characterization of finite minimal automorphic posets of width three}
\author{\sc Frank a Campo}
\date{\small Viersen, Germany\\
{\sf acampo.frank@gmail.com}}

\maketitle
\begin{abstract}
The characterization of the finite minimal automorphic posets of width three is still an open problem. Niederle has shown that this task can be reduced to the characterization of the nice sections of width three having a non-trivial tower of nice sections as retract. We solve this problem for a sub-class $\fN_2$ of the finite nice sections of width three. On the one hand, we characterize the posets in $\fN_2$ having a retract of width three being a non-trivial tower of nice sections, and on the other hand we characterize the posets in $\fN_2$ having a 4-crown stack as retract. The latter result yields a recursive approach for the determination of posets in $\fN_2$ having a 4-crown stack as retract. With this approach, we determine all posets in $\fN_2$ with height up to six having such a retract. For each integer $n \geq 2$, the class $\fN_2$ contains $2^{n-2}$ different isomorphism types of posets of height $n$.
\newline

\end{abstract}

\section{Introduction} \label{sec_Intro}

We call a finite poset $P$ {\em automorphic} if it has a fixed point free automorphism, and we call it {\em minimal automorphic} if additionally every proper retract of it has the fixed point property. Minimal automorphic posets have been present in the literature since the late eighties; an overview is provided by \cite{Schroeder_2022_MASoC}. In particular, they have found interest as forbidden retracts because a poset does not have the fixed point property iff it has a minimal automorphic poset as retract \cite[Th.\ 4.8]{Schroeder_2016}.

However, the description or characterization of minimal automorphic posets is still an open field. (Unless otherwise stated, ``poset'' always means ``finite poset'' in this article.) Since the pioneering work of Rival \cite{Rival_1976,Rival_1982}, we know that a poset of height one does not have the fixed point property iff it contains a crown, because a crown of minimal length will be a retract of it. Brualdi and Da Silva \cite{Brualdi_DdaSilva_1997} generalized this result by showing that a poset does not have the fixed point property iff it has a ``generalized crown'' as retract. With respect to posets of small width, Fofanova and Rutkowski \cite{Fofanova_Rutkowski_1987} showed that a poset of width two does not have the fixed point property iff it has a 4-crown stack as retract. (The result extends even to chain-complete infinite posets.) Niederle \cite{Niederle_1989} characterized posets of width three not having the fixed point property by having a non-trivial retract being a ``tower of nice sections''. Later on, he extended his result to posets of width four \cite{Niederle_2008}.

Unfortunately, it is quite difficult to break down Niederle's abstract characterization of the fixed point property to concrete posets or classes of posets of width three. According to our knowledge the only real success in this field is a result of Farley \cite{Farley_1997}, who showed that the ranked nice sections of width three are exactly the 6-crown stacks, and that a 6-crown stack is minimal automorphic iff its height $h$ is not a multiple of three. Otherwise, it has a 4-crown stack of height $\frac{2}{3} h$ as retract. Schr\"{o}der \cite[Lemma 2.6]{Schroeder_2022_MASoC} generalized one direction of this result by showing that a $2n$-crown stack of height $nm$ has a 4-crown stack of height $2m$ as retract. Moreover, in the case of $n \geq 3$, a $2n$-crown stack of height $h$ is not minimal automorphic if $h$ and $n$ have a non-trivial common factor \cite[Prop.\ 2.7]{Schroeder_2022_MASoC}.

In our article, we extend the characterization of minimal automorphic posets of width three to a class of posets which goes beyond the 6-crown stacks. For each integer $n \geq 2$, the class contains $2^{n-2}$ isomorphism types of posets of height $n$.

After the preparatory Section \ref{sec_preparation}, we recapitulate in Section \ref{sec_Sections} the definitions of {\em sections, nice sections,} and {\em towers of nice sections} introduced by Niederle \cite{Niederle_1989} for posets of width three. As mentioned above, he pointed out that the very heart of the characterization problem is the characterization of nice sections of width three having a non-trivial tower of nice sections as retract.

We start our investigation by characterizing the ``rather cumbersome to define nice sections'' \cite[p. 136]{Farley_1997} by means of standard terms in Proposition \ref{prop_irred}. Afterwards we introduce in Definition \ref{def_horizon} the {\em horizon} of a nice section $P$ of width three and height $h_P \geq 2$: With $P(0), \ldots , P(h_P)$ denoting the {\em level sets} of $P$, the {\em horizon} of $P$ is the smallest integer $\eta \in \myN$ for which, for all level-indices $0 \leq k < \ell \leq h_P$, the relation $\ell - k \geq \eta$ implies $P(k) < P(\ell)$, i.e., {\em all} points in level $P(k)$ are below {\em all} points in level $P(\ell)$. The horizon of a poset is thus the distance in which all details become blurred.

The subject of our investigation are nice sections with horizon two and height at least two, and we collect them in the class $\fN_2$. In Proposition \ref{prop_isomorphism} in Section \ref{sec_horizonTwo}, we characterize the isomorphism types of these posets. In the following Sections \ref{sec_RetrW3} and \ref{sec_retrSplit}, we deal with their retracts. Our approach is to look at {\em lower segments} and {\em upper segments} of $P \in \fN_2$ and to prove results of the type
\begin{quote}
$P$ has a retraction $r : P \rarr R$ with properties XYZ iff there exist a lower segment $L$ of $P$ and an upper segment $U$ of $P$ with retractions $ s : L \rarr S$ and $t : U \rarr T$ providing properties ABC and DEF, respectively.
\end{quote}
Using this approach, we characterize in Theorem \ref{theo_32} those posets $P \in \fN_2$ which have a retract of width three being a non-trivial tower of nice sections. Furthermore, we characterize in Theorem \ref{theo_split} the posets $P \in \fN_2$ having a 4-crown stack as retract. The latter result yields a recursive approach for the determination of posets in $\fN_2$ having a 4-crown stack as retract. With this approach, we determine in Section \ref{sec_Application} all posets in $\fN_2$ with height up to six having such a retract.

\section{Preparation} \label{sec_preparation}

In this section, we introduce our notation and recapitulate definitions of structures being in the focus of our investigation. For all other terms of order theory, the reader is referred to standard textbooks as \cite{Schroeder_2016}. 

Let $P = (X, \leq)$ be a finite poset. (We deal with finite posets only in this paper.) For $Y \subseteq X$, the {\em induced sub-poset} $P \vert_Y$ of $P$ is $\left( Y, {\leq} \cap (Y \times Y) \right)$. To simplify notation, we identify a subset $Y \subseteq X$ with the poset $P \vert_Y$ induced by it.

For $y \in X$, we define the {\em down-set} and {\em up-set induced by $y$} as
\begin{align*}
\darr_P y & := \mysetdescr{ x \in X }{ x \leq y }, \quad \odarr_P y := ( \darr_P y ) \setminus \setx{y}, 
\\
\uarr_P y & := \mysetdescr{ x \in X }{ y \leq x }, \quad \ouarr_P y := ( \uarr_P y ) \setminus \setx{y}. 
\end{align*}
If not required, we omit the labeling with $P$ in the notation of down-sets and up-sets. Given two points $x < y$ in $P$, the point $x$ is called a {\em lower cover} of $y$ and $y$ an {\em upper cover} of $x$ if $x \leq z \leq y$ implies $z \in \setx{x,y}$ for all $z \in P$; this relation is denoted by $x \lessdot y$. For integers $i \leq j$, we write $[i,j]$ for the interval $\setx{i, \ldots , j}$.

Let $A, B \subseteq P$. We write $A < B$ if $a < b$ for all $a \in A$, $b \in B$. If $A = \setx{a}$ is a singleton, we simply write $a < B$ and $B < a$. The notation for the relations $\leq, \lessdot, >, \geq, \gtrdot$ is analogous.

The {\em length} of a chain is its cardinality minus 1, and the height of a poset is the maximal length of a chain contained in it. For a poset $P$, we denote its height by $h_P$. The {\em width} of a poset is the largest size of an antichain contained in it.

For $n \in \myN \setminus \setx{1}$, a {\em $2n$-crown} is a poset with $2n$ points $x_{i,j}$, $i \in [0,n-1]$, $j \in \setx{0,1}$, in which $x_{0,0} < x_{0,1} > x_{1,0} < \cdots > x_{n-1,0} < x_{n-1,1} > x_{0,0}$ are the only comparabilities between different points. 

For a poset $P$, the {\em level sets} $P(k)$, $k \in [0,h_P]$, are recursively defined by
\begin{align*}
P(0) & := \min P, \\
P(k+1) & := \min \left( P \setminus \cup_{i=0}^k P(i) \right) \quad \mytext{for all } 0 \leq k < h_P.
\end{align*}
Additionally, we write for all $k, \ell \in [0,h_P]$ with $k \leq \ell$,
\begin{align*}
P(k,\ell) & := P(k) \cup P(\ell), \\
P(k \rightarrow \ell ) & := \cup_{i=k}^\ell P(i).
\end{align*}

Let $P$ and $Q_1, \ldots , Q_n$ be posets. We say that $P$ is a {\em $(Q_1, \ldots , Q_n)$-stack} if 
$$
h_P = \sum_{j=1}^n h_{Q_j} \quad \mytext{and} \quad 
P \left( \sum_{j=1}^{J-1} h_{Q_j} \rarr \sum_{j=1}^J h_{Q_j} \right) \simeq Q_J \; \; \mytext{for all } J \in [1, n].
$$
Moreover, for a poset $T$ with $h_T$ being a divisor of $h_P$, we say that $P$ is a {\em $T$-stack} if $P( i h_T \rightarrow (i+1) h_T) \simeq T$ for all $i \in [0, h_P / h_T - 1]$. All 4-crown stacks with the same height are isomorphic, and this holds also for all 6-crown stacks with the same height \cite[Prop.\ 4.1]{Farley_1997}.

Given pairwise disjoint posets $P_1 = (X_1, \leq_1), \ldots , P_n = (X_n, \leq_n)$, their {\em ordinal sum} $P_1\oplus \ldots \oplus P_n$ is defined as the poset having the carrier $\cup_{i=1}^n X_i$ and the partial order relation $\preceq$ defined as follows: for $x \in X_i$, $y \in X_j$, the relation $x \preceq y$ is equivalent to
$$
[ i = j \mytext{ and} x \leq_i y ] \; \; \mytext{ or } \; \; [i < j].
$$
The posets $P_i$ are called {\em ordinal summands} of $P_1\oplus \ldots \oplus P_n$. The level sets of an ordinal sum are exactly the level sets of the ordinal summands.

\begin{figure}
\begin{center}
\includegraphics[trim = 70 720 215 70, clip]{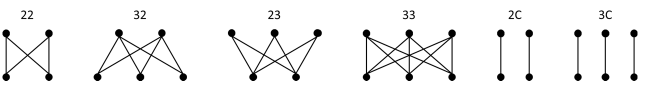}
\caption{\label{fig_Types_nm} Posets of type $22, 32, 23, 33, 2C$, and $3C$.}
\end{center}
\end{figure}

For $A$ and $B$ being antichains of size 2 or 3, we simply write ``$\#A \#B$'' as shortcut for the isomorphism type of $A \oplus B$, e.g., $Q \simeq 22$ is a 4-crown. Additionally, we write $2C$ and $3C$ for a poset consisting of two (three) disjoint chains of height 1 without any additional comparabilities between different points. The shortcuts are illustrated in Figure \ref{fig_Types_nm}.

For a mapping $f : X \rarr Y$, we denote by $f \vert_Z$ the pre-restriction of $f$ to a subset $Z$ of its domain and by $f \vert^Z$ its post-restriction to a subset $Z$ of its codomain with $f[X] \subseteq Z$. An order homomorphism $r : P \rightarrow P$ is called a {\em retraction} of the poset $P$ if $r$ is idempotent, and an induced sub-poset $R$ of $P$ is called a {\em retract of $P$} if a retraction $r : P \rightarrow P$ exists with $r[X] = R$. For the sake of simplicity, we identify $r$ with its post-restriction and write $r : P \rightarrow R$. A poset $P$ has the fixed point property iff every retract of $P$ has the fixed point property \cite[Th.\ 4.8]{Schroeder_2016}.

\begin{lemma} \label{lemma_retrOrdProd}
Let $Q = P_1 \oplus \ldots \oplus P_n$ be an ordinal sum and $r : Q \rightarrow R$ a retraction. $R$ has the fixed point property if $r[P_i] \not\subseteq P_i$ for an $i \in [1,n]$.
\end{lemma}
\BP Let $p \in P_i$ with $r(p) \in P_j$ for an index $j \not= i$. Assume $j > i$. For $x \in \cup_{k=0}^{j-1} P(k)$, we have $x < r(p)$, hence $r(x) \leq r(p)$, and for $x \in \cup_{k=j}^n P(k)$, we have $p < x$, thus $r(p) \leq r(x)$. For every endomorphism $f$ of $R$, we thus have $f(r(p)) \leq r(p)$ or $r(p) \leq f(r(p))$, and by the Theorem of Abian-Brown \cite{Abian_Brown_1961}, \cite[Th.\ 3.32]{Schroeder_2016}, $R$ has the fixed point property. The case $j < i$ is dual.

\EP

A point $a \in P$ is called {\em retractable} if $P \setminus \setx{a}$ is a retract of $P$. Retractability of $a$ is equivalent to the existence of a point $b \in P \setminus \setx{a}$ with $\odarr a \subseteq \darr b$ and $\ouarr a \subseteq \uarr b$, and $a$ is mapped to such a point by every retraction $r : P \rightarrow P \setminus \setx{a}$. Additionally, a point $a \in P$ is called {\em irreducible} if it has a single lower cover or a single upper cover. Irreducibility of a point is a special kind of retractability, and if the point $a \in P$ is irreducible, then $P$ has the fixed point property iff $P \setminus \setx{a}$ has the fixed point property \cite[Schol.\ 4.13]{Schroeder_2016}. 

A poset is called {\em automorphic} if it has a fixed point free automorphism, and it is called {\em minimal automorphic} if additionally every proper retract of it has the fixed point property. The following characterization of automorphic posets of width 3 goes back to \cite[Th.\ 10]{Niederle_2008}:

\begin{lemma}[{\cite[Prop.\ 4.38]{Schroeder_2016}}] \label{lemma_autom3} 
A poset $P$ of width at most three is automorphic iff for all level indices $0 \leq k < \ell \leq h_P$, the level-pair $P(k,\ell)$ is a 6-crown or of type 22, 23, 32, 33, $2C$, or $3C$.
\end{lemma}

\section{Sections and nice sections} \label{sec_Sections}

The following definition is due to Niederle \cite{Niederle_1989}:

\begin{definition} \label{def_Section}
A {\em section} is a two-element antichain or a poset $P$ of heigth at least one with carrier $X := \mysetdescr{ c_{k,j} } { k \in [0,h_P], j \in \nz }$ for which
\begin{itemize}
\item $C_j := c_{0,j} < \ldots < c_{h_P,j}$ is a chain for all $j \in \nz$;
\item $\setx{ c_{k,0}, c_{k,1}, c_{k,2} }$ is an antichain for all $k \in [0,h_P]$;
\item For all indices $k, \ell \in [0, h_P]$ and $i,j \in \nz$, we have $c_{k,i} < c_{\ell,j} \Rightarrow c_{k,i+1} < c_{\ell,j+1}$ (here as in what follows, index calculation is modulo 3);
\item $P(k,k+1) \not\simeq 33$ for all level indices $k \in [0, h_P-1]$.
\end{itemize}
A section $P$ is called {\em nice} if for all $x, y \in P$ with $x < y$
$$
\ouarr x \not\subseteq \uarr y \quad \mytext{and} \quad \odarr y \not\subseteq \darr x,
$$
and $\fN$ is the class of nice sections. A {\em tower of sections} is an ordinal sum of sections.
\end{definition}

For a section $P$ of width 3, the first and second property imply $c_{0,j} \lessdot \ldots \lessdot c_{h_P,j}$ for all $j \in \nz$. We call the chains $C_0, C_1, C_2$ the {\em main chains} of $P$. Clearly, $P(k) = \setx{ c_{k,0}, c_{k,1}, c_{k,2} }$ for all $k \in [0,h_P]$. We define mappings $\lambda : P \rightarrow [0,h_P]$ and $\gamma : P \rarr \setx{0,1,2}$ by setting $\lambda(c_{k,j}) := k$ and $\gamma(c_{k,j}) := j$, hence $\setx{p} = P(\lambda(p)) \cap C_{\gamma(p)}$ for all $p \in P$. Two consecutive levels of $P$ form a 6-crown or are of type $3C$. In what follows, the level-index function $\lambda$ and the chain-index function $\gamma$ refer always to the levels and main chains of $P$; for a retract $R$ of $P$, we thus have $\setx{x} = P(\lambda(x)) \cap C_{\gamma(x)}$ for all $x \in R$.

Examples for nice sections are the 6-crown stacks. Farley \cite{Farley_1997} showed that a 6-crown stack $P$ is minimal automorphic iff $h_P$ is not a multiple of three. Otherwise, it has a 4-crown stack of height $\frac{2}{3} h_P$ as retract.

We use the following two results of Niederle:

\begin{lemma}[{\cite[p.\ 121 and p.\ 125]{Niederle_1989}}] \label{lemma_Niederle}
The decomposition of a tower of sections into sections is unique, and a poset of width at most three has not the fixed point property iff it has a tower of nice sections as retract.
\end{lemma}

We conclude that an automorphic poset $P$ of width at most 3 not being a tower of nice sections is not minimal automorphic.

\begin{figure}
\begin{center}
\includegraphics[trim = 70 640 280 70, clip]{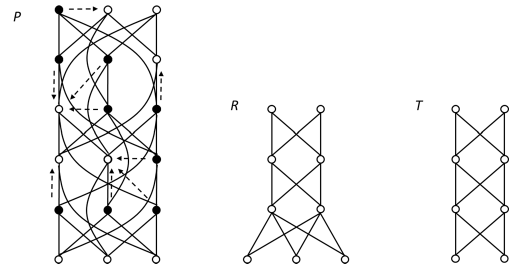}
\caption{\label{fig_Niederle} A nice section $P$ of width three with two retracts $R$ and $T$. The arrows indicate a retraction onto $R$. $R$ is not a tower of sections but $T$ is a tower of nice sections.}
\end{center}
\end{figure}

Figure \ref{fig_Niederle} shows an example of a nice section $P$ of width three and two retracts $R$ and $T$ of it. $R$ shows that $P$ is not minimal automorphic, but $R$ is not a tower of sections (the 3-antichain is the obstacle). However, $R$ has a 4-crown stack $T$ as retract, and this is a retract of $P$ being a tower of nice sections, as predicted by Lemma \ref{lemma_Niederle}. Niederle \cite{Niederle_1989} called sections ``very nice'' if they do not have a non-trivial tower of sections as retract and asked for their characterization. The problem is still open, too.

Let $P$ be a tower of nice sections $S_1, \ldots S_n$, thus $P = S_1 \oplus \cdots \oplus S_n$. In investigating if $P$ is minimal automorphic, we can due to Lemma \ref{lemma_retrOrdProd} restrict us to retractions $r : P \rightarrow R$ with $r \vert_{S_i}^{S_i} : S_i \rightarrow S_i$ being a retraction for all $i \in [1,n]$: $P$ is minimal automorphic iff $S_i$ is minimal automorphic for all $i \in [1,n]$. The objects we have to deal with in investigating minimal automorphic posets of width three are thus the nice sections of width three.

In the following proposition, we characterize the ``rather cumbersome to define nice sections'' \cite[p. 136]{Farley_1997} by standard terms:

\begin{proposition} \label{prop_irred}
Let $P$ be a section of width three. Equivalent are:
\begin{enumerate}
\item $P$ is nice;
\item $P$ does not contain a retractable point;
\item $P$ does not contain an irreducible point.
\end{enumerate}
\end{proposition}
\BP $1 \Rightarrow 2$: Let $a, b \in P$ be incomparable with $\lambda(a) \leq \lambda(b)$. In the case of $\lambda(a) < \lambda(b)$, we have $c_{\lambda(b), \gamma(a)} \in (\ouarr a) \setminus \uarr b$ and $c_{\lambda(a), \gamma(b)} \in (\odarr b) \setminus \darr a$.

Now let $\lambda(a) = \lambda(b)$. For $\lambda(a) < h_P$, $P(\lambda(a),\lambda(a)+1)$ is a 6-crown or of type $3C$, and $P(\lambda(a)+1)$ contains an upper cover of $a$ not being an upper cover of $b$ and vice versa. For $\lambda(a) > 0$, use the dual argument.

$2 \Rightarrow 3$ is trivial.

$3 \Rightarrow 1$: Let $a < b$ and let $x$ be an upper cover of $a$ with $a \lessdot x \leq b$. For the existing upper cover $y \not= x$ of $a$, we must have $b \not\leq y$. The point $b$ is treated dually.

\EP

\begin{figure}
\begin{center}
\includegraphics[trim = 70 690 370 70, clip]{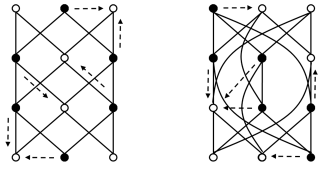}
\caption{\label{fig_N2h3} Two nice sections with 4-crown stacks as retracts.}
\end{center}
\end{figure}

Farley \cite[Prop.\ 4.1]{Farley_1997} showed that all 6-crown stacks of equal height are isomorphic. In their Hasse-diagrams, all 6-crowns can thus be drawn  in a standardized way with a ``void'' in the center as in Figure \ref{fig_N2h3} on the left. In Proposition \ref{prop_isomorphism} we will see that also in the Hasse-diagrams of the nice sections we are going to work with we can draw all 6-crowns in this way.

In a nice section $P$ not being a 2-antichain, the level-pairs $P(0,1)$ and $P(h_P-1,h_P)$ must form 6-crowns. The only nice sections of height 0, 1, and 2 are thus the 2-antichains, the 6-crowns, and the 6-crown stacks of height 2. All three are minimal automorphic. Two nice sections with height 3 are shown in Figure \ref{fig_N2h3}. Both contain a 4-crown stack of height 2 as retract, as indicated in the figure. (The retraction of the 6-crown stack is from \cite[Fig.\ 5.10]{Farley_1997}.)

In a section of width three, two consecutive level sets form a 6-crown or are of type 3C, whereas in a tower of sections, also 22, 23, 32, and 33 may appear. If two consecutive levels in a tower of {\em nice} sections are of type $23$, $32$, or $33$, we can look one level beyond the level set(s) being a 3-antichain:

\begin{corollary} \label{coro_2333_Bedeutung}
Let $P$ be a tower of nice sections with $P(k,k+1) \simeq$ 23 or 33 for a level index $k \in [0,h_P-1]$. The level sets $P(k)$ and $P(k+1)$ are level sets of different ordinal summands of $P$, and $P(k+1)$ is the bottom level of an ordinal summand of $P$. In particular, $k \leq  h_P-2$, and $P(k+1,k+2)$ is a 6-crown.
\end{corollary}
\BP The levels sets $P(k)$ and $P(k+1)$ are level sets of ordinal summands of $P$, and they cannot belong to the same nice ordinal summand due to the fourth item in Definition \ref{def_Section}. The 3-antichain $P(k+1)$ must be the bottom level of a nice section being an ordinal summand of $P$ and the result follows.

\EP

For $P \in \fN$ with height $h_P \geq 2$, we have $P(0,h_P) \simeq 33$. The following definition is thus meaningful:
\begin{definition} \label{def_horizon}
For a nice section $P \in \fN$ of height $h_P \geq 2$, we define the {\em horizon} of $P$ as the smallest integer $\eta\in \myN$ with
$$
P(k,k+\eta) \simeq 33 \quad \mytext{for all } k \in [0,h_P-\eta].
$$
For an integer $\eta \in \myN$, the symbol $\fN_\eta$ denotes the class of nice sections with height $\geq 2$ and horizon $\eta$.
\end{definition}

The horizon of a nice section is thus the distance in which all details become blurred. The set $\fN_1$ is empty, and the nice sections not contained in $\fN_\eta$ for any $\eta \geq 2$ are the 2-antichains and the 6-crowns. 

The horizon $\eta$ of a nice section $P$ has structural impact because for every $p \in P(k)$, all lower and upper covers of $p$ must belong to $P(k-\eta \rightarrow k-1)$ and $P(k+1 \rightarrow k+\eta)$, respectively. This restricts 3-element level sets of retracts being towers of nice sections:

\begin{lemma} \label{lemma_eta_Rell}
Let $P \in \fN_\eta$ with $\eta \geq 2$, let $R$ be a retract of $P$ being a tower of nice sections, and let $R(\ell) = \setx{x,y,z}$ be a 3-antichain in $R$ with $\lambda(x) \leq \lambda(y) \leq \lambda(z)$. Then $\lambda(z) - \lambda(x) \leq \eta-1$ with $\lambda(x) < \lambda(y) < \lambda(z)$ in the case of equality.
\end{lemma}
\BP The first inequality is a direct consequence of $\setx{x,z}$ being an antichain.

The level set $R(\ell)$ is a level set of a nice section $S$ of width 3 being an ordinal summand of $R$. Let $S(k) = R(\ell)$ and assume first that $k$ is not the height of $S$. Then $\# S(k+1) = 3$, and there exist $x', y', z' \in S(k+1)$ with $\gamma(x') = \gamma(x)$, $\gamma(y') = \gamma(y)$, and $\gamma(z') = \gamma(z)$.

Assuming $\lambda(z) - \lambda(x) = \eta - 1$, we get
\begin{align*}
\lambda(x) = \lambda(y) & \Rightarrow x, y, z < z', \\
\lambda(y) = \lambda(z) & \Rightarrow x < x', y', z',
\end{align*}
thus $S(k,k+1) = 33$ in both cases due to the third  property in Definition \ref{def_Section}. But this contradicts the fourth property in Definition \ref{def_Section}. For $k$ being the height of $S$, apply the dual argument on $S(k-1,k)$.

\EP

For later use, we notate a simple observation:

\begin{lemma} \label{lemma_hP_und_eta}
Let $\eta \geq 2$. A poset $P \in \fN_\eta$ with $h_P > \eta$ does not have a 4-crown as retract.
\end{lemma}
\BP Let $r : P \rightarrow R$ be a retraction to a 4-crown and let $a, b \in P(0)$, $v, w \in P(h_P)$ with $r[\setx{a,b}] = R(0)$, $r[\setx{v,w}] = R(1)$. A forbidden point $x \in P(1)$ exists with $\setx{a, b} \lessdot x < \setx{v, w}$.

\EP

Now we can easily see that all 4-crown stacks $R$ being a retract of one of the posets $P \in \fN_2$ in Figure \ref{fig_N2h3} have $R(0) \subset P(ß0)$ and $R(h_R) \subset P(3)$. (For the poset on the left, this is already stated in \cite[Lemma 5.9]{Farley_1997}.) Assume $\setx{a,b} = R(0) \not\subset P(0)$ with $\lambda(a) \leq \lambda(b)$. The sets $r^{-1}(a)$ and $r^{-1}(b)$ are disjoint and contain both at least one point of $P(0)$ which yields $a \in P(0)$ and $b \in P(1)$. But then $R(1 \rarr h_R) \subset P(2,3)$, and $R(1 \rarr h_R)$ cannot be a 4-crown stack. But $h_R = 1$ contradicts the last lemma.

\section{Nice sections with horizon two} \label{sec_horizonTwo}

For the rest of this article we will focus on nice sections with horizon two. In this section, we collect first results about them and their retracts.

\begin{proposition} \label{prop_isomorphism}
The isomorphism type of a poset $P \in \fN_2$ is uniquely determined by the types of its level-pairs $P(k, k+1)$, $k \in [1, h_P-2]$. In particular, there exist $2^{n-2}$ isomorphism types of posets of height $n \geq 2$ in $\fN_2$. Furthermore, every permutation of $P(0)$ can be uniquely extended to an automorphism of $P$.
\end{proposition}
\BP For each integer $n \geq 2$ and each index set $K \subseteq \myNkz{n-1}$ with $0, n-1 \in K$, we can construct a poset $P \in \fN_2$ with $h_P = n$ and $P(k, k+1)$ being a 6-crown iff $k \in K$: We start with the union $Q$ of three pairwise 
disjoint chains of length $n$, extend $Q(k,k+1)$ to an arbitrary 6-crown for every $k \in K$, and complete the partial order relation by adding $Q(k) \times Q(\ell)$ for all $k, \ell \in [0,n]$ with $k \leq \ell - 2$.

Now let $P = (X, \leq_P)$ and $Q = (Y,\leq_Q)$ be elements of $\fN_2$ with $h_P = h_Q$, and let $\alpha : P(0) \rarr Q(0)$ be a bijection. We define on $X$ and $Y$ the relations
\begin{align*}
\prec_P & := \mysetdescr{ (x,y) \in {<_P} }{ \lambda(y) - \lambda(x) = 1 } \\ 
\mytext{and} \quad \prec_Q & := \mysetdescr{ (x,y) \in {<_Q} }{ \upsilon(y) - \upsilon(x) = 1 },
\end{align*}
where $\upsilon$ is the level-function of $Q$. It is easily seen that we can extend $\alpha$ to an isomorphism $\beta : (X, \prec_P) \rarr (Y, \prec_Q)$ iff
\begin{equation*}
\forall \; k \in [1,h_P-2] : \quad 
P(k,k+1) \mytext{ 6-crown} \; \; \Leftrightarrow \; \; Q(k,k+1) \mytext{ 6-crown},
\end{equation*}
and that in this case $\beta$ is uniquely determined by $\alpha$. Every isomorphism from $(X, \prec_P)$ to $(Y, \prec_Q)$ sends $P(k)$ onto $Q(k)$ for all $k \in [0,h_P]$, and due to $P, Q \in \fN_2$, the isomorphisms between $(X, \prec_P)$ and $(Y, \prec_Q)$ are exactly the isomorphisms between $P$ and $Q$.

\EP

A direct consequence is

\begin{corollary} \label{coro_stacks}
Let $P \in \fN_2$, let $r : P \rarr R$ be a retraction with $r[P(0)] = R(0) \subseteq P(0)$, and let $f : P(0) \rarr P(0)$ be an idempotent mapping with $\# f[P(0)] = \# R(0)$. There exists a retraction $s : P \rarr S \simeq R$ with $s(x) = f(x)$ for all $x \in P(0)$.
\end{corollary}
\BP Let $R(0) = \setx{p}$ be a singleton and $f[P(0)] = \setx{q}$. According to Proposition \ref{prop_isomorphism}, there exists an automorphism $\sigma$ of $P$ with $\sigma(p) = q$, and $\sigma \circ r \circ \sigma^{-1}$ is a retraction with the desired property.

Now let $\# R(0) = 2$. Let $P(0) = \setx{a,b,z}$, $R(0) = \setx{a,b}$, $r(z) = b$, and let $f[P(0)] = \setx{a',b'}$, $z' \in P(0) \setminus \setx{a',b'}$, and $f(z') = b'$. There exists an automorphism $\sigma$ of $P$ with $\sigma(a) := a'$, $\sigma(b) := b'$, and $\sigma(z) := z'$, and again the mapping $\sigma \circ r \circ \sigma^{-1}$ has all desired properties.

For $\# R(0) = 3$, both mappings $f$ and $r$ are the identity mapping on $P(0)$.

\EP

In what follows, we frequently work with sub-posets $P(k_1 \rarr k_2)$ of posets $P \in \fN_2$. We call them {\em segments} in general, {\em lower segments} if $k_1 = 0$, and {\em upper segments} if $k_2 = h_P$. Due to Proposition \ref{prop_isomorphism}, all possible segments, lower segments, and upper segments are contained in the following classes:
\begin{align*}
\fNS_2 & := \mysetdescr{ P(1 \rarr h_P-1) }{ P \in \fN_2 }, \\
\fNL_2 & := \mysetdescr{ P(0 \rarr h_P-1) }{ P \in \fN_2 }, \\
\fNU_2 & := \mysetdescr{ P(1 \rarr h_P) }{ P \in \fN_2 }.
\end{align*}
We have $\fN_2 \subset \fNL_2 \cap \fNU_2$ and $\fNL_2 \cup \fNU_2 \subset \fNS_2$. Furthermore, Proposition \ref{prop_isomorphism} extends to $\fNS_2$: Two segments $U, V \in \fNS_2$ of equal height $n$ are isomorphic iff both are 3-antichains (the case $n = 0$) or all of their consecutive level-pairs are isomorphic (the case $n \geq 1$), and also the extension of a level permutation to an automorphism is always possible.

The most tedious part of our work will deal with retracts $R$ with $R(0)$ or $R(h_R)$ being a 2-antichain. We show in the Lemmata \ref{lemma_R0P0} and \ref{lemma_R0P0_C3} that without loss of generality we can assume that they and their retractions provide some standard features. We start with two simple observations; in particular the second one will play an important role in many proofs:

\begin{corollary} \label{coro_simpleResults}
Let $P \in \fNS_2$ and let $r : P \rarr R$ be a retraction.

1) For all $\ell \in [0,h_R]$ there exists an index $k \in [0,h_P-1]$ with $R(\ell) \subset P(k,k+1)$.


2) Let $\ell \in [0,h_R]$ with $\# R(\ell) =2$ and $1 \leq \rho := \max \mysetdescr{ \lambda(x) }{x \in R(\ell) } < h_P$. Then
\begin{align} \label{schubVonUnten}
\begin{split}
& r^{-1}(x) \cap P(0 \rarr \rho-1) \not= \emptyset \quad \mytext{for all points } x \in R(\ell) \\
\Rightarrow \quad & r[P(\rho+1 \rarr h_P)] = R(\ell+1 \rarr h_R) \subseteq P(\rho+1 \rarr h_P).
\end{split}
\end{align}
\end{corollary}
\BP The first part follows from the definition of the horizon. For the second part observe that the assumption implies $x \leq r(y)$ for all $x \in R(\ell)$ and all $y \in P(\rho+1 \rarr h_P)$, hence $R(\ell) < r(y)$. 

\EP

\begin{lemma} \label{lemma_R0P0}
Let $P \in \fNL_2$ with $h_P \geq 2$ and let $r : P \rarr R$ be a retraction with $R(0)$ being a 2-antichain. 
\begin{enumerate}
\item We have $R(0) \subset P(0,1)$ with $R(0) \cap P(0) \not= \emptyset$, and if $R(0) \cap P(1) \not= \emptyset$, then $h_R \geq 1$ and the poset $R( 1 \rarr h_R)$ is a retract of $P(2 \rarr h_P)$.
\item There exists a retraction $s : P \rarr S \simeq R$ with $S(0) \subset P(0)$, $s[P(0)] = S(0)$, and $s^{-1}(p) = \setx{p}$ for a point $p \in S(0)$.
\end{enumerate}
\end{lemma}
\BP 1. Let $\setx{a,b} = R(0)$ with $\lambda(a) \leq \lambda(b)$. $r^{-1}(a)$ and $r^{-1}(b)$ are disjoint down-sets in $P$, hence $a \in P(0)$. Now use both parts of Corollary \ref{coro_simpleResults}.

2. Let $\setx{a,b} = R(0) \not\subseteq P(0)$. Then $a \in P(0)$, $b \in P(1)$ due to the first part. There exists a point $b' \in P(0) \cap \darr b$, and we have $b' \not= a$ and $r(b') = b$. We define
\begin{align*}
s(x) := &
\begin{cases}
r(x), &\mytext{if } x \in P \setminus r^{-1}(b),\\
b', &\mytext{otherwise,}
\end{cases}
\end{align*}
This mapping is idempotent, and for $x < y$, only two cases have to be discussed:
\begin{itemize}
\item $r(x) = b \not= r(y)$: $s(x) = b' < b = r(x) \leq r(y) = s(y)$.
\item $r(x) \not= b = r(y)$: Is impossible due to $r(x) \leq r(y) = b \in R(0)$.
\end{itemize}
$s$ is thus a retraction with $S := s[P] \simeq R$ and $S(0) \subset P(0)$. Because $P(0,1)$ is a 6-crown, there must be a point $p \in S(0)$ with $s^{-1}(p) = \setx{p}$. In the case of $s[P(0)] \not= S(0)$ we redirect the point $z \in P(0) \setminus S(0)$ to the point in $S(0) \setminus \setx{p}$ and are finished.

\EP

\begin{figure}
\begin{center}
\includegraphics[trim = 70 715 375 70, clip]{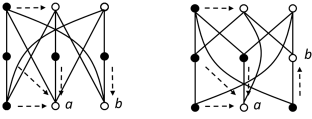}
\caption{\label{fig_P013C} Examples for the retractions described in Lemma \ref{lemma_R0P0_C3}.}
\end{center}
\end{figure}

\begin{lemma} \label{lemma_R0P0_C3}
Let $P \in \fNS_2$ with $h_P \geq 2$ and $P(0,1) \simeq 3C$ and let $r : P \rarr R$ be a retraction with $R(0)$ being a 2-antichain. 
\begin{enumerate}
\item $R(0) \subset P(0,1)$ and $r[P(2 \rarr h_P)] = R(1 \rarr h_S)$.
\item A retraction $s : P \rarr S \simeq R$ exists with $S(0) \subset P(0)$, $s^{-1}(S(0)) = P(0,1)$, and $s[P(2 \rarr h_P)] = S(1 \rarr h_S)$. (For an example, see the left part of Figure \ref{fig_P013C}.)
\item If $P(1,2)$ is a 6-crown, then a retraction $s : P \rarr S \simeq R$ exists with $S(0) \cap P(0) \not= \emptyset$, $S(0) \cap P(1) \not= \emptyset$, $s^{-1}(S(0)) = P(0,1)$, and $s[P(2 \rarr h_P)] = S(1 \rarr h_S)$. (An example is shown in the right part of Figure \ref{fig_P013C}.)
\end{enumerate}
\end{lemma}
\BP 1. Let $R(0) = \setx{a,b}$ with $\lambda(a) \leq \lambda(b)$. Because the sets $r^{-1}(a)$ and $r^{-1}(b)$ are disjoint down-sets in $P$, $\lambda(b) \geq 2$ is not possible. In the case of $\lambda(b) = 1$, the equation $r[P(2 \rarr h_P)] = R(1 \rarr h_P)$ follows with implication  \eqref{schubVonUnten}, and in the case of $\lambda(b) = 0$, the assumption $P(0,1) \simeq 3C$ yields $R(1) \subset P(2 \rarr h_P)$.

2. Let $a, b \in P(0)$ be two different points. Defining $s : P \rarr R$ by \begin{align} \label{retr_P01_C3}
s(x) & := 
\begin{cases}
r(x), & \mytext{if } x \in P(2 \rarr h_P), \\
a, & \mytext{if } x \in \setx{a, c_{1,\gamma(a)}}, \\
b, & \mytext{if } x \in P(0,1) \setminus \setx{a, c_{1,\gamma(a)}}
\end{cases}
\end{align}
yields according to the first part a retraction of $P$ with the desired properties.

3. Let $s$ be a retraction as described in the second part. There exists a point $p \in P(1)$ with $p < S(1)$. We define for all $x \in P$
$$
s'(x) :=
\begin{cases}
s(x), & \mytext{if } s(x) \not= s(p), \\
p, & \mytext{otherwise.}
\end{cases}
$$


\EP

\begin{corollary} \label{coro_P01_3C}
Let $P \in \fNS_2$ with $h_P \geq 2$ and $P(0,1) \simeq 3C$. P has a retract $R$ with $R(0,1)$ being a 4-crown iff $P(2 \rarr h_P)$ has a retract $T$ with $T(0)$ being a 2-antichain. 
\end{corollary}
\BP ``$\Rightarrow$'' is a direct consequence of Lemma \ref{lemma_R0P0_C3}. If $t : P(2 \rarr h_P) \rarr T$ is a retraction with $T(0)$ being a 2-antichain, select two points $a, b \in P(0)$ and define a retraction of $P$ as in \eqref{retr_P01_C3}.

\EP

\section{Non-trivial retracts of width three being towers of nice sections} \label{sec_RetrW3}

In this section, $P$ is a nice section belonging to $\fN_2$ and $r : P \rightarrow R$ is a retraction onto a retract $R$ of width three being a non-trivial tower of nice sections. In particular, no level set of $R$ is a singleton. As previously, the level-index function $\lambda$ refers always to the levels of $P$.

\begin{corollary} \label{coro_eta_subsec}
If $S$ is an ordinal summand of $R$ being a nice section of width 3, then $S = P(k \rightarrow \ell)$ for some $k < \ell$.
\end{corollary}
\BP Let $i \in [0,h_S]$ and $S(i) = \setx{x,y,z}$ with $\lambda(x) \leq \lambda(y) \leq \lambda(z)$. Lemma \ref{lemma_eta_Rell} yields $\lambda(z) - \lambda(x) \leq 1$, even $\lambda(z) = \lambda(x)$ because $\lambda(x) < \lambda(y) < \lambda(z)$ is not possible. There exist thus level-indices $k_0 < \ldots < k_{h_S}$ of $P$ with $S(i) = P(k_i)$ for all $i \in [0,h_S]$, and the horizon of $P$ enforces $k_{i+1} = k_i+1$ for all $i \in [0,h_S-1]$, because $S(i,i+1) \simeq 33$ is not possible (fourth item in Definition \ref{def_Section}).

\EP

\begin{lemma} \label{lemma_33}
$R(\ell,\ell+1) \not\simeq 33$ for all $\ell \in [0, h_R-1]$. In particular, $\# R(\ell) = 3$ for all $\ell \in [0,h_R]$ yields $R = P$.
\end{lemma}
\BP Let $R(\ell,\ell+1) \simeq 33$. We have seen in Corollary \ref{coro_2333_Bedeutung} and its dual that $R(\ell)$ is the top level set of a nice section $S$ of width 3 being an ordinal summand of $R$ and that $R(\ell+1)$ is the bottom level set of another nice section $T$  of width 3 being an ordinal summand of $R$. Corollary \ref{coro_eta_subsec} yields indices $k_1 < k_2 < k_3 < k_4$ with $S = P(k_1 \rightarrow k_2)$ and $T = P(k_3 \rightarrow k_4)$, and $33 \simeq R(\ell,\ell+1) = P(k_2, k_3)$ yields $k_3 \geq k_2 + 2$. Clearly, $r[ P(k_2 \rightarrow k_3) ] =  R(\ell,\ell+1)$.

Let $y \in P(k_2+1 \rarr k_3-1)$. We have $x < y$ for all $x \in P(k_2-1) = R(\ell-1)$, and because of $P(k_2-1,k_2) \not\simeq 33$, the point $r(y)$ cannot belong to $P(k_2) = R(\ell)$. In the same way we see that $r(y)$ cannot belong to $R(\ell+1)$, thus $r[P(k_2+1 \rarr k_3-1)] \cap R(\ell,\ell+1) = \emptyset$, a contradiction.

If $\# R(\ell) = 3$ for all $\ell \in [0,h_R]$, our result yields $R(\ell,\ell+1)$ being a 6-crown or of type $3C$ for all $\ell \in [0,h_R-1]$, and $R$ cannot be decomposed in an ordinal sum with two or more summands. $R$ is thus itself a nice section. Corollary \ref{coro_eta_subsec} delivers $i < j$ with $R = P(i \rightarrow j)$, and because $P(0,1)$ and $P(h_P-1,h_P)$ are 6-crowns we must have $i = 0$ and $j = h_P$.

\EP

The following lemma will cause a chain-reaction in the proof of Theorem \ref{theo_32}:

\begin{lemma} \label{lemma_schieben}
Assume that there exists a level index $\ell \in [0, h_R-1]$, with $R(\ell) < R(\ell+1)$. If there exists a level index $k \in [0, h_P-1]$ with
\begin{align*}
R(\ell) & = P(k) \\
\mytext{or} \quad R(\ell) & \subset P(k), \quad P(k,k+1) \simeq 3C, \quad \mytext{and} \quad r[P(k+1)] \subseteq R(\ell+1 \rarr h_R),
\end{align*}
then $k \leq h_P-3$ with
\begin{align*}
P(k+1,k+2) & \simeq 3C, \\
\mytext{and} \quad \quad \quad r[ P(k+1) ] & = R(\ell+1) \subseteq P(k+2).
\end{align*}
\end{lemma}
\BP We start with showing that under both assumptions $k \leq h_P-2$ and
\begin{equation} \label{eq_Enthaltenheiten}
r[P(k+1)] \; \; \subseteq \; \; R(\ell+1 \rarr h_R) \; \; \subseteq \; \; P(k+2 \rarr h_P).
\end{equation}

\begin{itemize}
\item $R(\ell) = P(k)$: According to Corollary \ref{coro_2333_Bedeutung}, the set $R(\ell)$ is the top-level of an ordinal summand $S$ of $R$ of width 3. Such a summand is a nice section, thus $S = P(j \rarr k)$ according to Corollary \ref{coro_eta_subsec}, and $R(\ell-1,\ell) = P(k-1,k)$ is a 6-crown. For every $x \in P(k+1)$, we have $R(\ell-1) < x$, but for every $p \in R(\ell)$, there exists a point $q \in R(\ell-1)$ with $q \not< p$. We conclude $r(x) \in R(\ell+1 \rarr h_R)$. Furthermore, $P(k) = R(\ell) < R(\ell+1)$ yields $k \leq h_P-2$ and $R(\ell+1 \rarr h_R) \subseteq P(k+2 \rarr h_P)$.
\item $R(\ell) \subset P(k)$, $P(k,k+1) \simeq 3C$, and $r[P(k+1)] \subseteq R(\ell+1 \rarr h_R)$: Again $R(\ell) < R(\ell+1)$ yields $k \leq h_P-2$ and $R(\ell+1 \rarr h_R) \subseteq P(k+2 \rarr h_P)$.
\end{itemize}

Now we can prove our result. Assume that $P(k+1,k+2)$ is a 6-crown or that $R(\ell+1)$ contains a point of $P(k+3 \rarr h_P)$. Then, on the one hand, $P(k+1) \cup R(\ell+1)$ is connected, and on the other hand, there exists for every $x \in P(k+1)$ a point $y \in R(\ell+1)$ with $x < y$. But then the first inclusion in \eqref{eq_Enthaltenheiten} yields $r[ P(k+1) \cup R(\ell+1)] = R(\ell+1)$, a contradiction. We thus have $R(\ell+1) \subseteq P(k+2)$ and $P(k+1,k+2) \simeq 3C$, hence $k \leq h_P - 3$.

It remains to show $R(\ell+1) = r[P(k+1)]$. ``$\subseteq$'' is due to \eqref{eq_Enthaltenheiten}, and implication \eqref{schubVonUnten} delivers $R(\ell+2 \rarr h_R) \subseteq P(k+3 \rarr h_P)$, thus $P(k+1) < R(\ell+2)$.

\EP

\begin{figure}
\begin{center}
\includegraphics[trim = 70 620 330 70, clip]{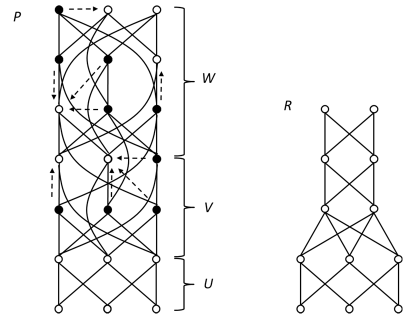}
\caption{\label{fig_PR_Width3} A poset $P \in \fN_2$ having a non-trivial retract $R$ of width three being a tower of nice sections. The segments $U$, $V$, and $W$ refer to Theorem \ref{theo_32}.}
\end{center}
\end{figure}

From Lemma \ref{lemma_33} we know that a non-trivial retract of $P \in \fN_2$ of width three being a tower of nice sections has to contain a level-pair of type $32$ or $23$. Nice sections with horizon two having such a retract are characterized in the following theorem and its dual. An example is shown in Figure \ref{fig_PR_Width3}.

\begin{theorem} \label{theo_32}
$P$ has a retract $R$ being a tower of nice sections with $R(\ell, \ell+1) \simeq 32$ for a level index $\ell$ of $R$, iff $P$ is a $(U,V,W)$-stack with \begin{itemize}
\item $U$ is a 6-crown or an element of $\fN_2$.
\item $V \in \fNS_2$ is a section with $h_V \geq 2$ even and $V(1 \rarr h_V)$ being a $3C$-stack.
\item $W \in \fN_2$ with $h_W \geq 3$ has a retraction $s : W \rightarrow S$ where $S$ is a tower of nice sections with $S(0,1) \simeq 22$, $s[W(1)] \cap S(0) = \emptyset$, and $h_S \geq 2$.
\end{itemize}
\end{theorem}
\BP ``$\Rightarrow$'': {\bf Existence of $U$:} $R(\ell)$ is a 3-antichain and belongs according to Corollary \ref{coro_2333_Bedeutung} to an ordinal summand $S$ of $R$ of width 3. Such a summand is a nice section, thus $S = P(j \rarr k)$ according to Corollary \ref{coro_eta_subsec}. We define $U := P(0 \rarr k)$.

{\bf Existence of $V$:} We apply Lemma \ref{lemma_schieben} iteratively.

{\em Initialization:} Due to $R(\ell) = P(k)$ and $R(\ell,\ell+1) \simeq 32$, Lemma \ref{lemma_schieben} delivers
\begin{align*}
P(k+1,k+2) & \simeq 3C, \\
r[ P(k+1) ] & = R(\ell+1) \subseteq P(k+2), \\
\mytext{and} \quad R(\ell+1) & \mytext{is a 2-antichain.}
\end{align*}

{\em Iteration:} Now assume that we have proven for an integer $J \in \myN$, that
\begin{align} \nonumber
P(k+1 \rarr k+2J ) & \mytext{is a } 3C\mytext{-stack,} \\ \label{32_3CStack}
r[ P(k+2j-1) ] & = R(\ell+j) \subseteq P(k+2j)  \quad \quad \mytext{for all } j \in [1,J], \\ \nonumber
R(\ell+1 \rarr \ell+J) & \mytext{is a 2-antichain or a 4-crown stack.}
\end{align}
Of course, $k+2J < h_P$, hence $\ell+J < h_R$.

In the case of $R(\ell+J,\ell+J+1) \simeq 23$, the Corollaries \ref{coro_2333_Bedeutung} and \ref{coro_eta_subsec} deliver an integer $I$ with $R(\ell+J+1) = P(k+2J+I)$. The dual of Lemma \ref{lemma_schieben} yields
$$
r[ P(k+2J+I-1)] = R(\ell+J) \subseteq P(k+2J+I-2).
$$
The comparison with \eqref{32_3CStack} delivers $I = 2$, but then we have
$$
r[P(k+2J+1)] = R(\ell+J) \stackrel{\eqref{32_3CStack}}{=} r[ P(k+2J-1)]
$$
in contradiction to that $P(k+2J-1,k+2J+1)$ is connected. Therefore,
\begin{equation} \label{RellJ_ellJ1_22}
R(\ell+J,\ell+J+1) \simeq 22.
\end{equation}
Furthermore, \eqref{32_3CStack} and implication \eqref{schubVonUnten} yield $r[P(k+2J+1)] \subseteq R(\ell+J+1 \rarr h_R)$. All together:
\begin{align*}
R(\ell+J) & < R(\ell+J+1), \\
R(\ell+J) & \subset P(k+2J), \\
\mytext{and} \quad \quad r[P(k+2J+1)] & \subseteq R(\ell+J+1 \rarr h_R).
\end{align*}
If now $P(k+2J,k+2J+1) \simeq 3C$, then Lemma \ref{lemma_schieben} yields
\begin{align*}
P(k+1 \rarr k+2J+2 ) & \mytext{is a } 3C\mytext{-stack,} \\
r[ P(k+2j-1) ] & = R(\ell+j) \subseteq P(k+2j)  \quad \mytext{for all } j \in [1,J+1], \\
R(\ell+1 \rarr \ell+J+1) & \mytext{is a 4-crown stack.}
\end{align*}

{\em Finalization:} The iteration eventually stops with an integer $K \in \myN$ for which $P(k+2K,k+2K+1)$ is a 6-crown. We define $V := P(k \rarr k+2K)$.

{\bf Existence of $W$:} Finally, we define $W := P(k+2K \rarr h_P) \in \fN_2$. Due to $r[P(k+2K-1)] \subset W(0)$ we have $r[W] \subset W$, and the mapping $s := r \vert_W^W$ is a retraction of $W$. The retract $S := s[W] = R(k+K \rarr h_R)$ is a tower of nice sections because $R$ is such a tower and $R(\ell+K-1, \ell+K) \simeq 32$ or $22$ for $K=1$ or $K>1$, respectively. Now $S(0)$ being a 2-antichain yields $h_W \geq 3$. Furthermore, implication \eqref{schubVonUnten} yields $s[W(1)] \cap S(0) = \emptyset$.   $S(0,1) = R(\ell+K,\ell+K+1) \simeq 22$ has already been seen in \eqref{RellJ_ellJ1_22}, and Lemma \ref{lemma_hP_und_eta} yields $h_S \geq 2$.

\begin{figure}
\begin{center}
\includegraphics[trim = 70 730 320 70, clip]{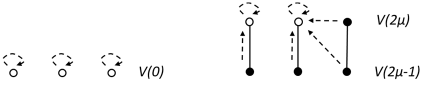}
\caption{\label{fig_V_Retr} The retraction $t$ of $V$ used in the proof of direction ``$\Leftarrow$'' of Theorem \ref{theo_32}. $1 \leq \mu \leq h_V / 2$.}
\end{center}
\end{figure}

``$\Leftarrow$'': Without loss of generality, assume $P(0 \rarr h_U) = U$, $P(h_U \rarr h_U + h_V) = V$, and $P(h_U + h_V \rarr h_P) = W$. Let $Y_j := y_{0,j} \lessdot \cdots \lessdot y_{h_V,j}$, $j \in \setx{0,1,2}$, be the main chains of $V$. We define a retraction $t : V \rarr T$ of $V$ by setting for all $(i,j) \in [0, h_V] \times \setx{0,1,2}$ (cf.\ Figure \ref{fig_V_Retr})
\begin{align*} 
t(y_{i,j}) & := 
\begin{cases}
y_{i,j}, & \mytext{if } i = 0, \\
y_{i+1,j}, & \mytext{if } i \mytext{ is odd and } j \in \setx{0,1}, \\
y_{i+1,1}, & \mytext{if } i \mytext{ is odd and } j = 2, \\
y_{i,j}, & \mytext{if } i > 0 \mytext{ is even and } j \in \setx{0,1}, \\
y_{i,1}, & \mytext{if } i > 0 \mytext{ is even and } j = 2, \end{cases}
\end{align*}
and we have $t \vert_{V(0)} = \id_{V(0)}$ and $t[V(h_V)] = T(h_T) \subset V(h_V)$.

Let $s$ be a retraction of $W$ as described in the assumptions. $s[W(1)] \cap S(0) = \emptyset$ implies $S(0) \subset W(0)$, and due to Lemma \ref{lemma_R0P0} and Corollary \ref{coro_stacks}, we can additionally assume $s[W(0)] = S(0)$ and $s(x) = t(x)$ for all $x \in P(h_U + h_V)$. Setting
\begin{align*}
r(x) & :=
\begin{cases}
x, & \mytext{if } x \in P(0 \rarr h_U), \\
t(x) & \mytext{if } x \in P(h_U \rarr h_U + h_V), \\
s(x), & \mytext{if } x \in P(h_U + h_V \rarr h_P)
\end{cases}
\end{align*}
yields a well-defined idempotent mapping.

Let $x, y \in P$ with $x < y$. The only case of interest is $x \in P(h_U + h_V - 1)$, $y \in P(h_U + h_V + 1)$. But then $t(x) \in T(h_T) = S(0)$ and $s(y) \in S(1 \rarr h_S)$, thus $r(x) < r(y)$ because $S(0,1)$ is a 4-crown. $r[ P(h_U \rarr h_U+2) ] \simeq 32$ is clear, and obviously $r[P]$ is a tower of nice sections.

\EP

The isomorphism type of the segment $V \in \fNS_2$ is uniquely determined by its height and the type of $V(0,1)$. Furthermore, $W$ cannot be a 6-crown stack because such a poset does not have a retraction as described in the theorem.

The description of the $(U,V,W)$-stack $P$ in Theorem \ref{theo_32} easily yields that $h_P \geq 6$ is required for a retract of width three being a non-trivial tower of nice sections. The poset in Figure \ref{fig_PR_Width3} is thus a smallest one in $\fN_2$ having such a retract. Three other isomorphism types of posets of height six in $\fN_2$ have such a retract, too. The first one we obtain by replacing in the poset in Figure \ref{fig_PR_Width3} the 6-crown $P(1,2)$ by a level-pair of type $3C$, the two other ones are the duals of these posets.

In the same way we see that a nice section $P \in \fN_2$ having a tower of nice sections as retract containing two different nice sections of width 3 as ordinal summands must have a height of at least 9, and even of at least 11 if it contains two different maximal 4-crown stacks as ordinal summands. Using Figure \ref{fig_PR_Width3} as a blueprint it is easy to construct examples.

\section{Four-crown stacks as retracts} \label{sec_retrSplit}

\begin{definition} \label{def_split}
Let $P \in \fN_2$. We call a quadrupel $(k, U, s, t)$ a {\em retractive up-split} of $P$ if
\begin{align*}
k & \in [0,h_P-1], \\
U &\subset P(0 \rarr k), \\
s & : P(0 \rarr k) \setminus U \rarr S \mytext{ is a retraction,} \\
\mytext{and} \quad \quad \quad \quad 
t & : P(k+1 \rarr h_P) \rarr T \mytext{ is a retraction} \\
\mytext{with } S \mytext{ and } T & \mytext{being 2-antichains or 4-crown stacks,}
\end{align*}
and we call a quadrupel $(k, D, s, t)$ a {\em retractive down-split} of $P$ if
\begin{align*}
k & \in [1,h_P], \\
D &\subset P(k \rarr h_P), \\
s & : P(k \rarr h_P) \setminus D \rarr S \mytext{ is a retraction,} \\
\mytext{and} \quad \quad \quad \quad 
t & : P(0 \rarr k-1) \rarr T \mytext{ is a retraction} \\
\mytext{with } S \mytext{ and } T & \mytext{being 2-antichains or 4-crown stacks.}
\end{align*}
\end{definition}
In order to avoid repetitions, we use the symbols $S$ and $T$ in what follows always as in this definition: for a retractive up-split $(k,U,s,t)$ of $P \in \fN_2$, we always have $S = s[P(0 \rarr k) \setminus U)]$ and $T = t[ P(k+1 \rarr h_P)]$, and correspondingly for a retractive down-split. The reader will observe that in a retractive up-split $(0,U,s,t)$ and a retractive down-split $(h_P,D,s,t)$ we must have $\# U \leq 1$ and $\# D \leq 1$, respectively, because $S$ contains at least two points.

Our main result in this section is
\begin{theorem} \label{theo_split}
A poset $P \in \fN_2$ has a 4-crown stack as retract iff a retractive up-split $(k,U,s,t)$ of $P$ exists with
\begin{align} \label{upsplit_begingungen}
\begin{split}
S(h_S) & <_P T(0),\\
\forall \; u \in U \;  \exists \; v \in T(0) :\; u & \not<_P p \mytext{ for all } p \in t^{-1}(v),
\end{split}
\end{align}
or a retractive down-split $(k,D,s,t)$ of $P$ exists with
\begin{align} \label{downsplit_begingungen}
\begin{split}
T(h_T) & <_P S(0),\\
\forall \; d \in D \;  \exists \; v \in T(h_T) :\; p & \not<_P d \mytext{ for all } p \in t^{-1}(v).
\end{split}
\end{align}
\end{theorem}

The reader will realize that the conditions \eqref{upsplit_begingungen} and \eqref{downsplit_begingungen} enforce $U \subseteq P(k)$ and $D \subseteq P(k)$, respectively. The key to the proof of the theorem is the following result:

\begin{proposition} \label{prop_RST}
Let $P \in \fN_2$ and let $(k,U,s,t)$ be a retractive up-split of $P$. There exists a retraction
\begin{align} \label{decomp_split_of_r}
\begin{split}
r : P & \rarr S \oplus T, \\
x & \mapsto
\begin{cases}
s(x), & \mytext{if } x \in P(0 \rarr k) \setminus U, \\
r(x) \in T, & \mytext{if } x \in U, \\
t(x), & \mytext{if } x \in P(k+1 \rarr h_P),
\end{cases}
\end{split}
\end{align}
iff the retractive up-split $(k,U,s,t)$ fulfills \eqref{upsplit_begingungen}.
\end{proposition}
\BP ``$\Rightarrow$'': $S \oplus T$ being an induced sub-poset of $P$ implies $S(h_S) <_P T(0)$. The second condition is a direct consequence of $r(u) \in T$ for all $u \in U$.

``$\Leftarrow$'':  Let $u \in U$. There exists a point $\tau(u) \in T(0)$ with $u \not<_P p$ for all $p \in t^{-1}(\tau(u))$. We define $\rho(u)$ as the single point contained in $T(0) \setminus \setx{ \tau(u) }$. Now we define a mapping $r : P \rarr P$ by setting for $x \in P$
\begin{align*}
r(x) & :=
\begin{cases}
s(x), & \mytext{if } x \in P(0 \rarr k) \setminus U, \\
\rho(x) , & \mytext{if } x \in U, \\
t(x), & \mytext{if } x \in P(k+1 \rarr h_P),
\end{cases}
\end{align*}
The mapping $r$ is clearly idempotent and because of $S(h_S) <_P T(0)$, the poset $S \oplus T$ is an induced sub-poset of $P$. It remains to show that $r$ is order preserving.

Let $x <_P y$. Only three cases have to be discussed:
\begin{itemize}
\item $x \in P(0 \rarr k) \setminus U$ and $y \in U \cup P(k+1 \rarr h_P)$ yields $r(x) \in S <_{S \oplus T} T \ni r(y)$.
\item $x \in U$, $y \in P(0 \rarr k)$ is not possible due to $U \subseteq P(k)$.
\item $x \in U$ and $y \in P(k+1 \rarr h_P)$: $x \not< p$ for all $p \in t^{-1}(\tau(x))$ delivers $t(y) \in T \setminx{\tau(x)}$. Because $\rho(x)$ is the only minimal point of this poset, we get $r(x) = \rho(x) \leq t(y) = r(y)$.
\end{itemize}

\EP

For the proof of Theorem \ref{theo_split}, it remains to show that a retraction $r' : P \rarr R'$ with $R'$ being a 4-crown stack guarantees the existence of a retraction $r : P \rarr R \simeq R'$ which can be decomposed by a retractive up-split as in \eqref{decomp_split_of_r} or by a retractive down-split as in the dual of \eqref{decomp_split_of_r}. We will even show more: every level set of $R'$ gives raise to a retractive up- or down-split of $P$. In order to avoid repetitions, $P$ is for the rest of this section an element of $\fN_2$ and $r : P \rarr R$ is a retraction onto a 4-crown stack.

\begin{corollary} \label{coro_Pkk}
If there exist level indices $\ell$ of $R$ and $k$ of $P$ with $R(\ell) \subset P(k)$, then at least one of the following retractive splits of $P$ exists:
\begin{itemize}
\item a retractive up-split $(\max \setx{k-1,0},U,s,t)$ fulfilling \eqref{upsplit_begingungen};
\item a retractive down-split $(\min \setx{k+1,h_P},D,s,t)$ fulfilling \eqref{downsplit_begingungen}.
\end{itemize}
\end{corollary}
\BP Let $z \in P(k) \setminus R(\ell)$. Three cases are possible:
\begin{enumerate}
\item $\ell \in [1, h_R-1]$: Then $k \in [1, h_P-1]$. In the case of $r(z) \in P(k \rarr h_P)$, we define $t := r \vert_{P(k \rarr h_P)}$,
$$
U := \mysetdescr{ x \in P(0 \rarr k-1)}{ r(x) \in P(k \rarr h_P)},
$$
$s := r \vert_{P(0 \rarr k-1)\setminus U}$, and Proposition \ref{prop_RST} yields the result. And in the case of $r(z) \in P(0 \rarr k-1)$, we construct a retractive down-split in the dual way.
\item $\ell = 0$: Then $k = 0$ due to the first part of Lemma \ref{lemma_R0P0}. If $r(z) \notin R(0)$, define
\begin{equation*}
t := r \vert_{P(1 \rarr h_P)}, \quad U := \setx{z}, \quad s := r \vert_{P(0)\setminus U}.
\end{equation*}
For $x \in P(1 \rarr h_P)$, we have $R(0) < x$ or $z < x$, hence $r(x) \in P(1 \rarr h_P)$. The mapping $t$ is thus a well defined retraction of $P(1 \rarr h_P)$ and Proposition \ref{prop_RST} yields the result.

Furthermore, in the case of $r(z) \in R(0)$, define
$$
t := r \vert_{P(0)}, \quad
D := \mysetdescr{ x \in P(1 \rarr h_P)}{ r(x) \in P(0)}, \quad
s := r \vert_{P(1 \rarr h_P)\setminus D},
$$
and apply the dual of Proposition \ref{prop_RST}.
\item $\ell = h_R$: Dual to the previous case.
\end{enumerate}

\EP

It remains to look at level sets $R(\ell)$ containing points from two consecutive level sets $P(k-1)$ and $P(k)$. According to the second part of Lemma \ref{lemma_R0P0} and its dual, we can assume without loss of generality $\ell \in [1, h_R-1]$ and $k \in [2, h_P-1]$. Now we can finish the proof of Theorem \ref{theo_split} in two steps:

\begin{lemma} \label{lemma_Pkk1}
Assume that there exists a level index $k \in [2, h_P-1]$ and points $d \in P(k-1)$, $e \in P(k)$ with $R(\ell) = \setx{d,e}$ for an $\ell \in [1, h_R-1]$. There exists a retraction $s: P \rarr R$ with 
\begin{equation*}
s^{-1}(d) \cap P(k,k+1) \not= \emptyset \quad \mytext{or} \quad s^{-1}(e) \cap P(k-1) \not= \emptyset.
\end{equation*}
\end{lemma}
\BP Let $x \in P(k-1)$ be a lower cover of $e$. In the case of $r(x) = e$ there is nothing to show. Assume thus $r(x) < d$ and let $U$ denote the set of upper covers of $x$. In the case of $r[U] \subseteq \uarr_R e$, the mapping $s$ with $s(p) := r(p)$ for $p \not= x$ and $s(x) := e$, yields a retraction because of $r[ \odarr_P x ] \subseteq P(0 \rarr k-2) < e$.

Let thus $u \in U$ with $r(u) \leq d$. In the case of $u \in P(k)$, the level-pair $P(k-1,k)$ must be a 6-crown due to $u \not= e$. But then $d \lessdot u$, hence $r(u) = d$. And in the case of $u \in P(k+1)$, the relation $d < u$ yields $r(u) = d$.

\EP

\begin{corollary} \label{coro_Pkk1}
Assume as in Lemma \ref{lemma_Pkk1} that there exists a level index $k \in [2, h_P-1]$ and points $d \in P(k-1)$, $e \in P(k)$ with $R(\ell) = \setx{d,e}$ for an $\ell \in [1, h_R-1]$. There exists a retractive up-split $(k,U,s,t)$ fulfilling \eqref{upsplit_begingungen} or a retractive down-split $(k-1,D,s,t)$ fulfilling \eqref{downsplit_begingungen}.
\end{corollary}
\BP According to Lemma \ref{lemma_Pkk1}, we can assume that at least one of the intersections $r^{-1}(d) \cap P(k,k+1)$ and $r^{-1}(e) \cap P(k-1)$ is non-empty.

In the case of $r^{-1}(e) \cap P(k-1) \not= \emptyset$, we have $r[P(k+1 \rarr h_P)] \subset P(k+1 \rarr h_P)$ due to implication \eqref{schubVonUnten}. Setting $t := r \vert_{P(k+1 \rarr h_P)}$,
$$
U := \mysetdescr{ x \in P(0 \rarr k)}{ r(x) \in P(k+1 \rarr h_P)},
$$
and $s := r \vert_{P(0 \rarr k) \setminus U}$, we get a retractive up-split $(k,U,s,t)$ fulfilling \eqref{upsplit_begingungen}. 

In the case of $r^{-1}(e) \cap P(k-1) = \emptyset$, we have $r^{-1}(d) \cap P(k,k+1) \not= \emptyset$ which implies $r[P(0 \rarr k-2)] \subset P(0 \rarr k-2)$ due to the dual of implication \eqref{schubVonUnten}. We construct a retractive down-split $(k-1,D,s,t)$ fulfilling \eqref{downsplit_begingungen} by setting $t := r \vert_{P(0 \rarr k-2)}$,
$$
D := \mysetdescr{ x \in P(k-1 \rarr h_P)}{ r(x) \in P(0 \rarr k-2)},
$$
and $s := r \vert_{P(k-1 \rarr h_P) \setminus D}$.

\EP

Now Theorem \ref{theo_split} is proven. As consequence we get:

\begin{figure}
\begin{center}
\includegraphics[trim = 70 705 300 70, clip]{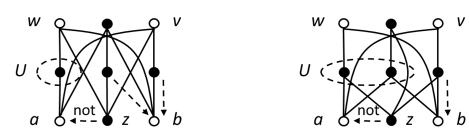}
\caption{\label{fig_CoroLuecke} Illustrations for the proof of Corollary \ref{coro_stackMitLuecke}. The three levels show $P(k \rarr k+2)$ for choices $b \in P(k)$, $w \in P(k+2)$ (which are not mandatory). The point $z$ is not mapped to $a$ by the retraction $s$ of $P(0 \rarr k)$.}
\end{center}
\end{figure}

\begin{corollary} \label{coro_stackMitLuecke}
Let $k \in [0,h_P-2]$ and assume that we have retractions $s : P(0 \rarr k) \rarr S$ and $t : P(k+2 \rarr h_P) \rarr T$ with $S$ and $T$ being 2-antichains or 4-crown stacks with
\begin{align*}
s^{-1}(a) & = \setx{a} \mytext{ for a point } a \in S(h_P), \\
\mytext{and} \quad t^{-1}(v) & = \setx{v} \mytext{ for a point } v \in T(0).
\end{align*}
The 4-crown stack $S \oplus T$ is a retract of $P$ if $P(k \rarr k+2)$ is not a 6-crown stack.
\end{corollary}

\BP Let $S(h_S) = \setx{a,b}$ and $T(0) = \setx{v,w}$. The assumptions imply $a \in P(k)$ and $v \in P(k+2)$. Assume $P(k+1,k+2) \simeq 3C$ and let $P(k) = \setx{a,b,z}$.

Let $U \subset P(k+1)$ be the set of upper covers of $a$ in $P(k+1)$. Applying Proposition \ref{prop_isomorphism} on $P(0 \rarr k+1)$, we see that without loss of generality we can assume that none of the points in $U$ is below $v$ (Figure \ref{fig_CoroLuecke}).

No point of $P(k+1) \setminus U$ is above $a$ and the point $b$ is the only maximal point of $S \setminx{a}$. Due to $s^{-1}(a) = \setx{a}$, the mapping
\begin{align*}
s' : P(0 \rarr k+1) \setminus U & \rarr S, \\
x & \mapsto 
\begin{cases}
s(x), & \mytext{if } x \in P(0 \rarr k), \\
b, & \mytext{if } x \in P(k+1) \setminus U
\end{cases}
\end{align*}
is a retraction, and due to $t^{-1}(v) = \setx{v}$, the quadrupel $(k+1, U, s', t)$ is a retractive up-split fulfilling \eqref{upsplit_begingungen}.

In the case of $P(k,k+1) \simeq 3C$, apply the dual construction.

\EP

\section{Application} \label{sec_Application}

In Theorem \ref{theo_split}, we deal with a {\em $t$-base} $V$ and an {\em $s$-base} $W$ of a poset $P \in \fN_2$, i.e.:
\begin{itemize}
\item $V$ is a lower or upper segment of $P$ having a retraction $t$ onto a 2-antichain or a 4-crown stack.
\item $W$ is an upper segment of $P$ having a retraction of $W \setminus D$ onto a 2-antichain or a 4-crown stack for a subset $D \subseteq W(0)$, or $W$ is a lower segment of $P$ having a retraction of $W \setminus U$ onto a 2-antichain or a 4-crown stack for a subset $U \subseteq W(h_W)$.
\item $V$ and $W$ are coupled by $V \cap W = \emptyset$ and $P = V \cup W$.
\end{itemize}

Assume that we have listed all lower segments up to height $n$ in $\fNL_2$ which have a 4-crown stack as retract. Given a poset $P \in \fN_2$ of height $n+1$, we can quickly identify all candidates for $t$-bases in $P$ by checking which of the lower segments in our list are lower segments $P(0 \rarr k)$ of $P$ and which of their duals are upper segments $P(k \rarr n+1)$ of $P$. Of course, $P(0)$ and $P(n+1)$ are always candidates for a $t$-base. According to Theorem \ref{theo_split}, the poset $P$ has a 4-crown stack as retract iff one of the corresponding segments $P(k+1 \rarr n+1)$ or $P(0 \rarr k-1)$ provides an $s$-base fulfilling \eqref{downsplit_begingungen} or \eqref{upsplit_begingungen}, respectively.

Using this recursive approach, we will identify all posets $P \in \fNL_2$ with height up to six having a 4-crown stack as retract. For some of them, e.g. the 6-crown stacks, we already know the result. We include them in our investigation because our purpose is to test our approach with as many nice sections as possible.

We encode the lower segments by the sequence of their level-pair types (cf.\ Proposition \ref{prop_isomorphism}). 1 indicates a 6-crown, 0 a type $3C$. The lower segment 111 is thus the 6-crown stack of height 3, whereas 110 indicates the lower segment starting with two 6-crowns followed by a single $3C$.

\begin{figure}
\begin{center}
\includegraphics[trim = 70 680 190 70, clip]{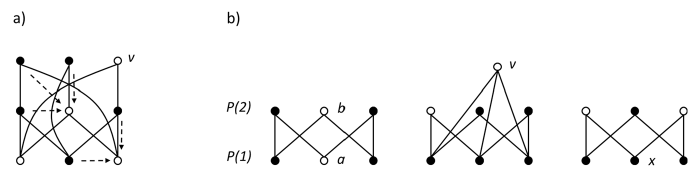}
\caption{\label{fig_Retr10Crit5} a) The lower segment $10$ and a retraction $r$ to a 4-crown with $r^{-1}(v) = \setx{v}$ for a point $v$ from the top level of the crown. b) Illustrations for the case discrimination in the proof of Criterion \ref{crit_split}.5. The hollow dots belong to $S(0)$.}
\end{center}
\end{figure}

\begin{table}
\begin{center}
{\footnotesize
\begin{tabular}{| l l c | l l c | l l c | l l c | }
\hline
1 & & n & & & & & & & & &  \\
\hline
11 & & n & 10 & & y & & & & & & \\
\hline
111 & & y & 101 & & y & 110 & & n & 100 & & n \\
\hline
1111 & 0,3 & n & 1101 & 0 & n & 1011 & 0,2,3 & n & 1001 & 0,2 & y \\
1110 & & n & 1100 & & n & 1010 & & y & 1000 & & y \\
\hline
11111 & 0,3 & n & 11101 & 0,3 & y & 11011 & 0 & n & 10111 & 0,2,3 &  y \\
11001 & 0 & y & 10101 & 0,2,3,4 & y & 10011 & 0,2,4 & y & 10001 & 0,2,4 & y \\
11110 & & y & 11100 & & y & 11010 & & n & 10110 & & y \\
11000 & & n & 10100 & & y & 10010 & & n & 10000 & & n \\
\hline
111111 & 0,3 & y & 111101 & 0,3,5 & y & 111011 & 0,3,5 & n & 110111 & 0 & n \\
101111 & 0,2,3,5 & y & 111001 & 0,3,5 & n & 110101 & 0 & n & 101101 & 0,2,3,5 & y \\
110011 & 0,5 & n & 101011 & 0,2,3,4,5 & n & 100111 & 0,2,4,5 & n & 110001 & 0 & n \\
101001 & 0,2,3,4,5 & y & 100101 & 0,2,4 & y & 100011 & 0,2,4,5 & n & 100001 & 0,2,4 & y \\
111110 & & n & 111100 & & n & 111010 & & n & 110110 & & n \\
101110 & & n & 111000 & & n & 110100 & & n & 101100 & & n \\
110010 & & n & 101010 & & y & 100110 & & y & 110000 & & n \\
101000 & & y & 100100 & & y & 100010 & & y & 100000 & & y \\
\hline
\end{tabular} }
\caption{\label{table_LS6} The lower segments from $\fNL_2$ with height up to six. The letters ``y'' and ``n'' indicate if a segment has a 4-crown stack as retract or not. Additional explanation in text.}
\end{center}
\end{table}

The lower segments up to height six are listed in Table  \ref{table_LS6}. The letters ``y'' and ``n'' indicate the result of our investigation whether the segment has a 4-crown stack as retract or not. For the posets 1, 10, 11, 111, and 101, we refer to common knowledge. The lower segment $10$ even has a retraction $r$ to a 4-crown with $r^{-1}(v) = \setx{v}$ for a point $v$ from the top level of the crown, as shown in Figure \ref{fig_Retr10Crit5}a.

For the posets with a final 0, we can decide by applying the dual of Corollary \ref{coro_P01_3C} referring to previous results in Table \ref{table_LS6}. As an example, the posets $110$ and $100$ do not have a 4-crown stack as retract, because the 6-crown $P(0,1)$ does not have a suitable retract. And $1010$ and $1000$ have a 4-crown stack as retract because 10 is marked with ``y''.

The bulk of the work concerns thus the lower segments $P$ with a 6-crown as final level-pair. They all are elements of $\fN_2$. For them, the integers in Table \ref{table_LS6} indicate the levels $k \in [0, h_P-1]$ for which $P(0 \rarr k)$ can be a $t$-base according to the results for lower segments with less height. For $1111$ it is 0 and 3: 0, because the antichain $P(0)$ is always a candidate for a $t$-base, and 3, because the poset $1111(0 \rarr 3) = 111$ is marked with ``y'' in the previous row. Because 1 and 11 are marked with ``n'' , the level indices 1 and 2 are missing in the list for $1111$.

The following criteria will be useful in our investigation:

\begin{criteria} \label{crit_split}
Let $P \in \fN_2$ with $h_P \geq 3$.
\begin{enumerate}
\item For every retractive down-split $(h_P,D,s,t)$ we have $t[P(0 \rarr h_P-3)] = T(0 \rarr h_T-1)$. In particular, the choice of $P(h_P)$ as an $s$-base implies that the poset $P(0 \rarr h_P-3)$ has a 2-antichain or a 4-crown tower as retract.
\item For the choice of $P(h_P-1,h_P)$ as $s$-base, the set $D \subseteq P(h_P-1)$ must contain at least two points.
\item The level-pair $P(h_P-1,h_P)$ is not a candidate for an $s$-base if $P(h_P-2,h_P-1)$ is a 6-crown.
\item $P(1 \rarr h_P)$ can serve as an $s$-base only together with a subset $D \subset P(1)$ containing at most a single point.
\item If $P(1,2)$ is a 6-crown, then $P(1 \rarr h_P)$ being an $s$-base implies that $P(3 \rarr h_P)$ has a 2-antichain or a 4-crown stack as retract.
\end{enumerate}
\end{criteria}
\BP 1. The first condition in \eqref{downsplit_begingungen} yields $T(h_T) \not\subset P(h_P-1)$. Now the dual of the first parts of the Lemmata \ref{lemma_R0P0} and \ref{lemma_R0P0_C3} yields $t[P(0 \rarr h_P-3)] = T(0 \rarr h_T - 1)$.

2. The poset $P(h_P-1,h_P) \setminus D$ has to be disconnected.

3. Is implied by the second criterion, because if $P(h_P-2,h_P-1)$ is a 6-crown, the third condition in \eqref{downsplit_begingungen} can be fulfilled for at most a single point $d \in P(h_P-1)$.

4. Trivial.

5. We discriminate three cases (cf.\ Figure \ref{fig_Retr10Crit5}b). If $S(0)$ contains a point $a \in P(1)$, then it has to contain a point $b \in P(2)$, too, because of the first condition in \eqref{downsplit_begingungen}. At least one of the points in $P(1) \setminus \setx{a}$ has to be mapped to $b$ (fourth criterion), thus $s[P(3 \rarr h_P)] = S(1 \rarr h_S)$ due to implication \eqref{schubVonUnten}.

Now assume that $S(0)$ does not contain a point of $P(1)$. If $S(0)$ contains a point $v \in P(3 \rarr h_P)$, then at least two points in $P(1)$ are below $S(0)$ in contradiction to the fourth criterion. And in the case of $S(0) \subset P(2)$, exactly one point $x \in P(1)$ is below $S(0)$ and has to be sent to $P(0)$. From the two remaining points in $P(1)$ each is under a single different point of $S(0)$. Due to the fourth criterion, they have to be mapped onto $S(0)$, and $s[P(3 \rarr h_P)] = S(1 \rarr h_S)$ follows again with implication \eqref{schubVonUnten}.

\EP

\begin{figure}
\begin{center}
\includegraphics[trim = 70 670 240 70, clip]{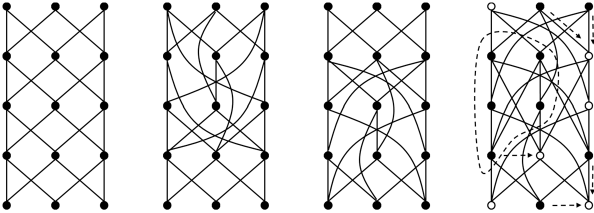}
\caption{\label{fig_PosetsHeight4} The posets $1111$, $1101$, $1011$, and $1001$. The poset $1011$ is also treated in \cite[p.\ 99]{Schroeder_2016} and mentioned in \cite[p.\ 136]{Farley_1997}. For the poset $1001$, a retraction onto a 4-crown tower is indicated resulting from a retractive down-split $(3,D,s,t)$ based on the retraction $t$ of $P(0 \rarr 2) = 10$ shown in Figure \ref{fig_Retr10Crit5}a.}
\end{center}
\end{figure}

The rightmost poset in Figure \ref{fig_PosetsHeight4} confirms that we can neither drop the condition ``$P(h_P-2,h_P-1)$ is a 6-crown'' in Criterion \ref{crit_split}.3 nor the condition ``$P(1,2)$ is a 6-crown'' in Criterion \ref{crit_split}.5.

\begin{figure}
\begin{center}
\includegraphics[trim = 70 660 190 70, clip]{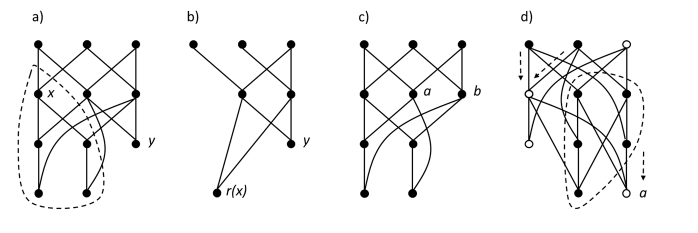}
\caption{\label{fig_011_001} The first three diagrams are illustrations for the proof of Lemma \ref{lemma_Q011}. Additionally the rightmost diagram shows for $P = 001$ a retraction of $P \setminus D$ to a 4-crown with a singleton $D \subset P(0)$.}
\end{center}
\end{figure}

\begin{lemma} \label{lemma_Q011}
Let $P = 011$. If $r : P \setminus D \rarr R$ is a retraction onto a 4-crown stack or a 2-antichain for a down-set $D$, then $D$ contains at least two points of $P(0)$.
\end{lemma}

\BP Assume that a retraction exists for $D = \setx{d}$ with $d \in P(0)$. Figure \ref{fig_011_001}a shows $P \setminx{d}$. Because $r$ cannot be extended to $P$, both upper covers $x$ and $y$ of $d$ have to be mapped to different points of $R(0)$. The encircled set $\darr x$ has thus to be mapped to a single point in $R(0)$. The resulting poset is shown in Figure \ref{fig_011_001}b. Clearly, it does not have a 4-crown stack or a 2-antichain as retract.

Now assume $D = \setx{d,y}$ with $d \in P(0)$ and $d \lessdot y$. The poset $P \setminus D$ is shown in Figure \ref{fig_011_001}c. If a retraction as described in the lemma exists, it cannot be extended to $P \setminx{d}$. The upper covers $a$ and $b$ of $y$ are thus mapped to different points of $R(0)$. But that is impossible because of $P(0) \setminx{d} < \setx{a,b}$.

\EP

Now we can start with the posets of height four. They are shown in Figure \ref{fig_PosetsHeight4}.

\begin{itemize}
\item P = 1111: Criterion \ref{crit_split}.1 prevents $P(4)$ from being a candidate for an $s$-base because the 6-crown $P(0,1)$ does not have a suitable retract. Furthermore, selecting $P(1 \rarr 4)$ as $s$-base does not work due to Criterion \ref{crit_split}.5, because $P(3,4) \simeq 1$ is marked with ``n'' in Table \ref{table_LS6}. $1111$ does thus have no retractive down-split fulfilling \eqref{downsplit_begingungen}, and because the poset is self-dual, we conclude that it does not have a 4-crown stack as retract.
\item $P = 1101$: Again, $P(1 \rarr 4)$ fails as an $s$-base due to Criterion \ref{crit_split}.5. $1101$ thus does not have a retractive down-split fulfilling \eqref{downsplit_begingungen}.
\item $P = 1011$: $P(1 \rarr 4)$, $P(3,4)$, or $P(4)$ being a successful $s$-base is prevented by Criterion \ref{crit_split}.4 together with Lemma \ref{lemma_Q011}, Criterion \ref{crit_split}.3, and Criterion \ref{crit_split}.1, respectively. Also $1011$ does thus not have a retractive down-split fulfilling \eqref{downsplit_begingungen}. Therefore, neither $1101$ nor $1011$ has a 4-crown stack as retract.


\item $P = 1001$: A retraction to a 4-crown stack resulting from a retractive down-split $(3,D,s,t)$ is indicated in Figure \ref{fig_PosetsHeight4} on the right.
\end{itemize}

We continue with the posets of height five:

\begin{itemize}
\item $P = 11111$: Due to the Criteria \ref{crit_split}.5 and \ref{crit_split}.3, the upper segments $P(1 \rarr 5)$ and $P(4,5)$ cannot work as $s$-bases. The self-dual poset has thus no 4-crown stack as retract.
\item $P = 11101$: see $P = 10111$.
\item $P = 11011$: Criterion \ref{crit_split}.5 prevents $P(1 \rarr 5)$ from being a successful $s$-base. Because $11011$ is self-dual, the poset cannot have a 4-crown stack as retract.
\item $P = 10111$: Figure \ref{fig_N2h3} shows that the 6-crown stack $Q = 111$ has a retraction $t : Q \rarr T$ with $T$ being a 4-crown stack and $t^{-1}(a) = \setx{a}$ for a point $a \in T(0)$. With $s : P(0) \rarr P(0)$ being a mapping with $s[P(0)]$ being a 2-antichain, apply Corollary \ref{coro_stackMitLuecke} with $k = 0$.
\item $P = 11001$: see $10011$.
\item $P = 10101, 10001$: $P(0 \rarr 2) = 10$ and $P(3 \rarr 5) = 01$ have retractions onto 4-crowns as shown in Figure \ref{fig_Retr10Crit5}a and its dual. By Corollary \ref{coro_stackMitLuecke} with $k = 2$, they can be combined to a retraction of $P$ to a 4-crown stack of height three.
\item $P = 10011$: Let $Q := 1001$. The retraction $t $ of $Q$ indicated in Figure \ref{fig_PosetsHeight4} on the right and a self-mapping $s$ of $P(5)$ onto a 2-antichain can be combined to a retractive down-split $(5, \emptyset, s, t)$ of $P$ fulfilling \eqref{downsplit_begingungen}.
\end{itemize}

For some of the posets of height six, a quick positive decision is possible:

\begin{itemize}
\item $P = 111111, 111101, 101111, 101101$: Combine the retractions in Figure \ref{fig_N2h3}.
\item $P = 101001, 100001$: Let $P$ be one of these posets. With $Q := P(3 \rarr 6) = 001$, select a point $d \in Q(0)$ and let $s$ be the retraction of $Q \setminus \setx{d}$ onto a 4-crown $S$ shown in Figure \ref{fig_011_001} on the right with $a \in S \cap Q(0)$ as in the figure. Furthermore, let $t$ be the retraction of $P(0 \rarr 2) = 10$ onto a 4-crown $T$ as shown in Figure \ref{fig_Retr10Crit5}a with $v \in T \cap P(2)$ as in the figure. Applying Proposition \ref{prop_isomorphism} on $P(0 \rarr 2)$, we can assume $v \lessdot a$ and $v \not< d$ and get a retractive down-split $(3,\setx{d},s,t)$ fulfilling \eqref{downsplit_begingungen}.
\item $P = 100101$: dual to $101001$.
\end{itemize}



The remaining posets of height six are $111011$, $111001$, $101011$, $110001$ and their duals and the self-dual poset $110011$. They do not have a 4-crown stack as retract:

\begin{lemma} \label{lemma_P36_gehtNicht}
None of the posets $111011$, $111001$, $101011$, $110001$, and $110011$ has a 4-crown stack as retract.
\end{lemma}
\BP Firstly, let $P \in \fN_2$ be any nice section of height six and $R$ a retract of $P$ being a 4-crown stack. Every level set of $R$ has to be contained in a single level set or in two consecutive level sets of $P$, and we conclude that at least one of the posets $P(0 \rarr 3) \cap R$ and $P(3 \rarr 6) \cap R$ is a 4-crown stack.

Now let $P$ be any of the five posets and $r$ a retraction onto a 4-crown stack $R$. With $S := P(3 \rarr 6) \cap R$, we show in a first step that $S$ cannot be a 4-crown stack.

$P(0 \rarr 3) = 111, 110$: $S$ being a 4-crown stack yields a point $x \in P(3)$ with $v := r(x) \in P(0 \rarr 2) = 11$. Let $\setx{v,w}$ be the level set of $R$ containing $v$. $x \in P(3)$ yields $v \geq r[P(0,1)]$, hence $w \notin P(0,1)$. But then, $r(z) < \setx{v,w}$ holds for at least five of the points $z \in P(0,1)$, and $P(0,1)$ cannot be mapped onto a 2-antichain.

$P = 101011$: $S$ being a 4-crown stack requires due to Lemma \ref{lemma_Q011} two points $x, y \in P(3)$ to be mapped into $P(0 \rarr 2)$. Three cases are possible, each of it leading to a contradiction:
\begin{itemize}
\item $\setx{r(x),r(y)}$ is a level set of $R$: Let $z \in P(2)$ be the common lower cover of $x$ and $y$. Due to $P(0) < z$, the point $r(z)$ is the only minimal point of $R$.
\item $r(x) < r(y)$: Then $r(x) \in R(0)$, and $P(0) < x$ yields that $r(x)$ is the only minimal point of $R$.
\item $v := r(x) = r(y)$: Let $\setx{v,w}$ be the level set of $R$ containing $v$. Due to $r[P(0 \rarr 2)] \leq v$, the point $w$ must belong to $P(3 \rarr h_P)$, and due to $\lambda(v) \leq 2$, we must have $w \in P(3)$ and $v \in P(2)$. But then $w$ belongs to a level set of $R$ contained in $S$.
\end{itemize}

$S$ is thus not a 4-crown stack. This requires a level set $R(\ell) = \setx{a,b}$ with $a \in P(2)$ and $b \in P(3)$. Clearly, $\ell = 1$. Due to $P(1) < b$, no point of $P(1)$ can be mapped to $a$. And if a point of $P(1)$ is mapped to $b$, implication \eqref{schubVonUnten} yields $R(2 \rarr h_R)$ being a retract of $P(4 \rarr 6) = 11$ which is impossible. But if no point of $P(1)$ can be mapped to $R(1)$, then $r[P(0,1)] = R(0)$ which is impossible, too.

\EP


\end{document}